\documentclass[12pt]{amsart}
\usepackage{enumerate, amscd, amssymb, fullpage}

\newtheorem{theorem}{Theorem}[section]
\newtheorem{lemma}[theorem]{Lemma}
\newtheorem{proposition}[theorem]{Proposition}
\newtheorem{corollary}[theorem]{Corollary}
\theoremstyle{definition}
\newtheorem{definition}[theorem]{Definition}

\newtheorem{example}[theorem]{Example}

\newcommand{\LL}{\mathcal{L}}

\newcommand{\FF}{\mathcal{ F}}

\def\bR{{\mathbb R}}
\def\bH{{\mathbb H}}
\def\bE{{\mathbb E}}
\def\bS{{\mathbb S}}

\newcommand\Curv{\operatorname{R}}
\newcommand\Ric{\operatorname{Ric}}

\newcommand\tr{\operatorname{Tr}}

\newcommand\h{\operatorname{h}}
 
\newenvironment{noindenum}
{\begin{list}{\labelitemi}{\leftmargin=0em \itemindent=0em}}
{\end{list}}

\newcommand\Sec{\operatorname{K}}
\newcommand\rk{\operatorname{rk}}

\newcommand\Ker{\operatorname{Ker}}

\catcode`\@=\active 

\begin{document}

\title{On the curvature of metric contact pairs}
\author{Gianluca~Bande}
\address{Dipartimento di Matematica e Informatica, Universit{\`a} degli studi di Cagliari, Via Ospedale 72, 09124 Cagliari, Italy}
\email{gbande{\char'100}unica.it}
\author{David E. Blair}
\address{Department of Mathematics, Michigan State University, East Lansing, MI 48824--1027, USA}
\email{blair{\char'100}math.msu.edu}
\author{Amine~Hadjar}
\address{Laboratoire de Math{\'e}matiques, Informatique et
Applications, Universit{\'e} de Haute Alsace - 4, Rue de
Fr{\`e}res Lumi{\`e}re, 68093 Mulhouse Cedex, France}
\email{mohamed.hadjar{\char'100}uha.fr}

\thanks{The first author was supported by the Project \emph{Start-up giovani ricercatori}--$2009$ of Universit\`a degli Studi di Cagliari and by a Visiting Professor fellowship at the Universit\'e de Haute Alsace in June 2010 and in June 2011. The second and the third authors were supported by a Visiting Professor fellowship at the Universit\`a degli Studi di Cagliari in April 2011 and January 2010 respectively, financed by Regione Autonoma della Sardegna.}
\date{\today; MSC 2010 classification: primary 53C25; secondary 53B20, 53C12, 53B35}

\begin{abstract}
We consider manifolds endowed with metric contact pairs
for which the two characteristic foliations are orthogonal.
We give
some properties of the curvature tensor and in particular a formula for the Ricci curvature in the direction of the sum of the two Reeb vector fields.
This shows that metrics associated to normal contact pairs cannot be flat.
Therefore flat non-K\"ahler Vaisman manifolds do not exist.
Furthermore we give a local classification of metric contact pair manifolds whose curvature vanishes on the vertical subbundle. As a corollary we have that flat associated metrics can only exist if the leaves of the characteristic foliations are at most three-dimensional.

\end{abstract}

\maketitle

\section{Introduction}
A \textit{contact pair} on a smooth even-dimensional manifold $M$ is a pair of one-forms
$\alpha_1$ and $\alpha_2$ of constant and complementary classes,
for which $\alpha_1$ restricted to the leaves of
the characteristic foliation of $\alpha_2$ is a contact form and vice versa \cite{Bande1,BH}. The Reeb vector fields on these contact leaves determine two global vector fields $Z_1$ and $Z_2$ called the Reeb vector fields of the pair.
This notion was first introduced by Blair, Ludden and Yano~\cite{Blair2} under the name \emph{bicontact} in the context of Hermitian geometry, and further studied by Abe \cite{Abe}.

In \cite{BH2, BH4} the first and the third authors
constructed metrics  adapted to contact pairs as in metric contact geometry.
More precisely, a \emph{metric contact pair} on an even dimensional manifold
is a triple $(\alpha_1 , \alpha_2 , g)$, where $(\alpha_1 , \alpha_2)$ is a contact pair
with Reeb vector fields $Z_1$, $Z_2$, and $g$
 is an associated metric, i.e.
 a Riemannian metric such that $g(X, Z_i)=\alpha_i(X)$, for $i=1,2$, and for which the endomorphism field $\phi$, uniquely defined by $g(X, \phi Y)= (d \alpha_1 + d
\alpha_2) (X,Y)$, satisfies
\begin{equation*}
\phi^2=-Id + \alpha_1 \otimes Z_1 + \alpha_2 \otimes Z_2 \, .
\end{equation*}

Contact pairs always admit associated metrics for which the two characteristic foliations are orthogonal \cite{BH2} or, equivalently, whose structure tensor $\phi$ is decomposable (i.e. $\phi$ preserves the characteristic distributions of $\alpha_1$ and $\alpha_2$).

In this paper we prove the following classification theorem which is analogous to that of the second author \cite{Blair4} concerning metric contact manifolds with curvature vanishing on the vertical subbundle:

\textbf{Main Theorem.}
\textit{Let $M$ be a $(2h+2k+2)$-dimensional manifold endowed with a metric contact pair $(\alpha_1, \alpha_2, \phi, g)$ of type $(h,k)$ (with $h\geq 1$) and decomposable $\phi$. If the curvature $\Curv$ of the metric $g$ satisfies $\Curv _{X Y}Z_i=0$ ($i=1,2$),
then $M$ is locally isometric to $\bE^{h+1}\times \bS^h(4)\times \bE^{k+1}\times \bS^k(4)$ if $k\geq 1$ or $\bE^{h+1}\times \bS^h(4)\times \bE^{1}$ if $k=0$.}

If the manifold is complete, then its Riemannian universal covering is globally isometric to $\bE^{h+1}\times \bS^h(4)\times \bE^{k+1}\times \bS^k(4)$ if $k\geq 1$ or $\bE^{h+1}\times \bS^h(4)\times \bE^{1}$ if $k=0$.  In this statement we will understand that when $h$ (or $k$) is equal to 1, the $\bS^h(4)$ factor will just contribute another line to the Euclidean factor.

As a corollary we obtain that the only manifolds which can carry flat metric contact pairs are either four or six-dimensional
with metric contact pairs of type $(1,0)$ or $(1,1)$ respectively.
We also prove several formulae concerning
 the curvature tensor and the Ricci curvature of a metric associated to a contact pair
 with decomposable $\phi$.
In particular, we show that, on a $2n$-dimensional manifold endowed with such a
 structure, the Ricci curvature of the associated metric in the direction
 of the vector field $Z=Z_1 + Z_2$ is $n-1-\frac{1}{2} \tr \h^2$, where $\h= \frac{1}{2} \LL_Z \phi$
 and $\LL_Z$ is the Lie derivative along $Z$.

 An immediate consequence is the non-existence of flat metrics associated to normal contact pairs with decomposable endomorphism.
 This implies that the metric of a non-K\"ahler Vaisman structure on a smooth manifold cannot be flat. This is interesting
 since
 the property is local; until now this result was well known only for \emph{closed} manifolds
 (see \cite{Vaisman} and \cite[Proposition 2.5]{Ornea}).

 In the sequel we denote by $\Gamma(B)$ the space of sections of a vector bundle $B$, by $\tr$ the trace of an endomorphism field,
 and by $\nabla$ the Levi-Civita connection of a given metric.
All the differential objects considered are assumed to be smooth.

\section{Preliminaries on metric contact pairs}\label{s:prelim}

In this section we gather the notions concerning contact pairs
that will be needed in the sequel. We refer the reader to
\cite{Bande2, BGK, BH, BH2, BH3, BH4, BK} for further informations and several
examples of such structures.

\subsection{Contact pairs}\label{s:prelimcp}
A pair $(\alpha_1, \alpha_2)$ of 1-forms on a manifold is said
to be a \emph{contact pair} of type $(h,k)$ if (see \cite{Bande1,BH}):
\begin{eqnarray*}
&\alpha_1\wedge (d\alpha_1)^{h}\wedge\alpha_2\wedge
(d\alpha_2)^{k} \;\text{is a volume form},\\
&(d\alpha_1)^{h+1}=0 \; \text{and} \;(d\alpha_2)^{k+1}=0.
\end{eqnarray*}
Since the form $\alpha_1$ (respectively $\alpha_2$) has constant class
$2h+1$ (respectively $2k+1$), the distribution $\Ker \alpha_1 \cap \Ker
d\alpha_1$ (respectively
 $\Ker \alpha_2 \cap \Ker d\alpha_2$) is completely integrable and then it determines the so-called \emph{characteristic
  foliation} $\mathcal{F}_1$ (respectively $\mathcal{F}_2$) whose leaves are endowed with a contact form induced by $\alpha_2$ (respectively $\alpha_1$).

The equations
\begin{eqnarray*}
&\alpha_1 (Z_1)=\alpha_2 (Z_2)=1  , \; \; \alpha_1 (Z_2)=\alpha_2
(Z_1)=0 \, , \\
&i_{Z_1} d\alpha_1 =i_{Z_1} d\alpha_2 =i_{Z_2}d\alpha_1=i_{Z_2}
d\alpha_2=0 \, ,
\end{eqnarray*}
where $i_X$ is the contraction with the vector field $X$, determine uniquely the two vector fields $Z_1$ and $Z_2$, called \emph{Reeb vector fields}. Since they commute \cite{Bande1, BH}, they give rise to a locally free $\mathbb{R}^2$-action, called  the \emph{Reeb
action}.

The tangent bundle of a manifold $M$ endowed with a contact pair
can be split in different ways. For $i=1,2$, let $T\mathcal F _i$
be the subbundle determined by the characteristic foliation of
$\alpha_i$, $T\mathcal G_i$ the subbundle of $TM$ whose fibers are
given by $\ker d\alpha_i \cap \ker \alpha_1 \cap \ker \alpha_2$
and $\mathbb{R} Z_1, \mathbb{R} Z_2$ the line bundles determined
by the Reeb vector fields. Then we have the following splittings:
$$
TM=T\mathcal F _1 \oplus T\mathcal F _2 =T\mathcal G_1 \oplus
T\mathcal G_2 \oplus \mathcal V  ,
$$
where $\mathcal V =\mathbb{R} Z_1 \oplus \mathbb{R} Z_2$.
Moreover we have $T\mathcal F _1=T\mathcal G_1 \oplus \mathbb{R}
Z_2 $ and $T\mathcal F _2=T\mathcal G_2 \oplus \mathbb{R} Z_1 $.
\begin{definition}
We say that a vector field is \textit{vertical} if it is a section of $\mathcal V$ and \textit{horizontal} if it is a section of $T\mathcal G_1 \oplus
T\mathcal G_2$. The subbundles $\mathcal V$ and $T\mathcal G_1 \oplus T\mathcal G_2$ will be called \textit{vertical} and \textit{horizontal} respectively.
\end{definition}
Notice that $d\alpha_1$ (respectively $d\alpha_2$) is symplectic on the vector bundle $T\mathcal G_2$ (respectively $T\mathcal G_1$).

\begin{example}
Take $(\mathbb{R}^{2h+2k+2},\alpha_1 , \alpha_2)$ where $\alpha_1$~, $\alpha_2$ are the Darboux contact forms on $\mathbb{R}^{2h+1}$
and $\mathbb{R}^{2k+1}$ respectively.
\end{example}
This is also a local model for all contact pairs of type $(h,k)$. Hence a contact pair manifold is locally the product of two contact manifolds \cite{Bande1, BH}.

\subsection{Contact pair structures}
We recall now the notion of contact pair structure studied in \cite{BH2, BH3, BH4}.
\begin{definition}[\cite{BH2}]
A \emph{contact pair structure} on a manifold $M$ is a triple
$(\alpha_1 , \alpha_2 , \phi)$, where $(\alpha_1 , \alpha_2)$ is a
contact pair and $\phi$ a tensor field of type $(1,1)$ such that:
\begin{equation}\label{d:cpstructure}
\phi^2=-Id + \alpha_1 \otimes Z_1 + \alpha_2 \otimes Z_2 , \;
\phi Z_1=\phi Z_2=0
\end{equation}
where $Z_1$ and $Z_2$ are the Reeb vector fields of $(\alpha_1 ,
\alpha_2)$.
\end{definition}
One can see that $\alpha_i \circ \phi =0$ for $i=1,2$ and that the rank of $\phi$ is
equal to $\dim M -2$ .
Since we are also interested in the induced structures, we recall the
following:
\begin{definition}[\cite{BH2}]
The endomorphism $\phi$ is said to be \emph{decomposable} if
$\phi (T\mathcal{F}_i) \subset T\mathcal{F}_i$, for $i=1,2$.
\end{definition}
The condition for $\phi$ to be decomposable is equivalent to $\phi(T\mathcal{G}_i)= T\mathcal{G}_i$ for $i=1,2$.

If $\phi$ is decomposable, then $(\alpha_1 , Z_1 ,\phi)$ (respectively
$(\alpha_2 , Z_2 ,\phi)$) induces, on every leaf of $\mathcal{F}_2$ (respectively $\mathcal{F}_1$), an almost contact structure (see e.g. \cite{Blairbook}) consisting of a contact form, its Reeb vector field and a structure tensor, the restriction
 of $\phi$ to the leaf.

On a manifold $M$ endowed with a contact pair, there always exists
a decomposable endomorphism field $\phi$ satisfying \eqref{d:cpstructure}
(see \cite{BH2}).

As a trivial example one can take two contact manifolds $M_i$, $i=1,2$ with structure tensors $( \alpha_i , \phi_i)$, and consider the contact pair structure
$(\alpha_1 , \alpha_2 , \phi_1 \oplus \phi_2)$ on $M_1 \times M_2$. In \cite{BH3} we gave examples of contact pair structures with decomposable endomorphism which are not locally products.

In what follows, on a manifold $M$ endowed with a contact pair structure $(\alpha_1, \alpha_2, \phi)$, we will consider the tensor fields defined by:
\begin{alignat*}{1}
N^1 (X,Y)= & [\phi , \phi ](X, Y) +2 d\alpha_1 (X,Y) Z_1 +2 d\alpha_2 (X,Y) Z_2 ,\\
N^2 _i (X,Y) = & (\LL_{\phi X} \alpha_i) (Y) - (\LL_{\phi Y}\alpha_i)(X) ,  \; i=1,2 , \\
\h=&\frac{1}{2}  \LL_{Z} \phi  \, ,
\end{alignat*}
for all $X, Y \in \Gamma (TM)$, where
$Z=Z_1+Z_2$ and $[\phi, \phi]$ is the Nijenhuis tensor of $\phi$. The vanishing of $N^1$ gives exactly
the \emph{normality} of the pair \cite{BH3}, that is the integrability of the two almost complex structures $\phi \pm (\alpha_1 \otimes Z_2 - \alpha_2 \otimes Z_1)$. In this case, by \cite[Equation (3.5)]{BH3} we have the following:

\begin{proposition}
\label{prop:normal-lziphi-vanish}
If a contact pair structure $(\alpha_1, \alpha_2, \phi)$ with Reeb vector fields $Z_1$ and $Z_2$ is normal,
we have $N^2_1=N^2_2=0$ , $\LL_{Z_1}\phi=\LL_{Z_2}\phi=0$ and then $\h=0$ .
\end{proposition}

\subsection{Metric contact pairs}
On manifolds endowed with contact pair structures it is natural
to consider the following metrics:
\begin{definition}[\cite{BH2}]
Let $(\alpha_1 , \alpha_2 ,\phi )$ be a contact pair structure on
a manifold $M$, with Reeb vector fields $Z_1$ and $Z_2$. A
Riemannian metric $g$ on $M$ is said to be:
\begin{enumerate}[{i)}]
\item \emph{compatible} if $g(\phi X,\phi Y)=g(X,Y)-\alpha_1 (X)
\alpha_1 (Y)-\alpha_2 (X) \alpha_2 (Y)$ for all $X,Y \in \Gamma
(TM)$,
\item \emph{associated} if $g(X, \phi Y)= (d \alpha_1 + d
\alpha_2) (X,Y)$ and $g(X, Z_i)=\alpha_i(X)$, for $i=1,2$ and for
all $X,Y \in \Gamma (TM)$. \label{ass-metric}
\end{enumerate}
\end{definition}
An associated metric is compatible, but the converse is not true.
\begin{definition}[\cite{BH2}]
A \emph{metric contact pair (MCP)} on a manifold $M$ is a
four-tuple $(\alpha_1, \alpha_2, \phi, g)$ where $(\alpha_1,
\alpha_2, \phi)$ is a contact pair structure and $g$ an associated
metric with respect to it. The manifold $M$ will be called an \emph{MCP manifold} or simply an \emph{MCP}.
\end{definition}

For an MCP $(\alpha_1, \alpha_2, \phi, g)$ the endomorphism field $\phi$ is decomposable
if and only if the characteristic foliations $\mathcal{F}_1 ,
\mathcal{F}_2$ are orthogonal \cite{BH2}. In this case $(\alpha_i, \phi, g)$ induces a metric contact structure
on the leaves of $\mathcal{F}_j$ , for $j\neq i$ .

Using a standard polarization on the symplectic vector bundles $T\mathcal G_i$, one can see that for a given contact pair $(\alpha_1, \alpha_2)$ there always exist a decomposable endomorphism field $\phi$ and a metric $g$ such
that $(\alpha_1, \alpha_2, \phi, g)$ is an MCP (see \cite {BH2}).  Moreover we have:
\begin{proposition}\label{p:Ni=0}
Let $(\alpha_1, \alpha_2, \phi , g)$ be an MCP
 with decomposable $\phi$. Then we have:
\begin{equation}\label{eq-lemma-N1=0}
N^2 _1=N^2 _2=0 .
\end{equation}
\end{proposition}
\begin{proof}
Since $\phi$ is decomposable, if $X, Y$ are vector fields tangent to different foliations, we have $d\alpha_i (\phi X, Y)=d\alpha_i (\phi Y, X)=0$, $i=1,2$ . If $X, Y$ are tangent to the same foliation $\FF_i$~, because $\phi$ preserves the foliation, we have $N_i ^2(X,Y)=0$. Moreover, for $j \neq i$, the triple $(\alpha_j , \phi , g)$ restricted to the leaves of $\FF_i$ is a metric contact structure and then it satisfies \eqref{eq-lemma-N1=0}, which is a well known fact in metric contact geometry. \end{proof}

Some other properties of MCP's are given by the
following results:
\begin{theorem}[\cite{BH2}]\label{th:ass-metrics}
Let $M$ be a manifold endowed with a contact pair structure $(\alpha_1 , \alpha_2, \phi)$, with Reeb vector fields $Z_1 , Z_2$. Let $g$ be
 a metric compatible metric with  the structure.
 Then we have:
\begin{enumerate}
\item $g(Z_i, X)= \alpha_i (X)$ for $i=1,2$ and for every $X \in \Gamma(TM)$;

\item $g(Z_i, Z_j)= \delta_{ij}$ for $i,j=1,2$;

\item $\nabla _{Z_i} Z_j = 0$ for $i,j=1,2$ (in particular the integral curves of the
Reeb vector fields are geodesics);

\item the Reeb action is totally geodesic (i.e. the orbits are totally geodesic two-dimensional submanifolds).

\end{enumerate}
Moreover, if $g$ is an associated metric, then $\LL_{Z_i}\phi=0$ if
and only if $Z_i$ is Killing.
\end{theorem}

In the normal case, by Proposition \ref{prop:normal-lziphi-vanish}, an immediate consequence is:

\begin{corollary}\label{cor:normal-killing}
If an MCP $(\alpha_1, \alpha_2, \phi, g)$ is normal, the Reeb vector fields $Z_1$ and $Z_2$ are Killing.
\end{corollary}

Now using the invariance of the $\alpha_i$'s under the flow of $Z=Z_1+Z_2$ one can also prove the following:
\begin{proposition}\label{p:Z-killing}
Let  $(\alpha_1 , \alpha_2, \phi, g)$ be an MCP with Reeb vector fields $Z_1$ and $Z_2$. Then $\h$ vanishes if and only if $Z$ is Killing.
\end{proposition}
We end this section with a result from \cite{BH4}:
\begin{theorem}[\cite{BH4}]
On an MCP manifold $(M, \alpha_1, \alpha_2, \phi, g)$ with decomposable $\phi$ the leaves of the characteristic foliations of the contact pair are orthogonal and minimal.
\end{theorem}
As example, one can simply take the product of two metric contact manifolds.
Here is an interesting example from \cite{BH4} which shows that an MCP manifold is not always locally the product of two metric contact manifolds:
\begin{example}\label{ex:nilp-lie-group}
Let us consider the simply connected $6$-dimensional nilpotent Lie group $G$ with structure
equations:
\begin{eqnarray*}
&d\omega_3= d\omega_6=0 \; \; , \; \; d\omega_2= \omega_5 \wedge
\omega _6 ,\\
&d\omega_1=\omega_3 \wedge \omega_4 \; \; , \; \;
d\omega_4= \omega_3 \wedge \omega_5 \; \; , \; \; d\omega_5 =
\omega_3 \wedge \omega_6 \, ,
\end{eqnarray*}
where the $\omega_i$'s form a basis for the cotangent space of $G$
at the identity.
Then $(\omega_1 , \omega_2)$ together with
the metric
\begin{equation*}
g=\omega_1 ^2+\omega_2 ^2+\frac{1}{2}\sum_{i=3}^6 \omega_i ^2
\end{equation*}
is a left invariant MCP of type $(1,1)$ on $G$. Note that the two characteristic foliations are orthogonal, and
that their leaves, although minimal, are not totally geodesic. So the metric $g$ is not locally a product.
Since the structure constants of the group are rational, there exist lattices $\Gamma$ such that $G/\Gamma$ is
compact. This MCP descends to all quotients $G/\Gamma$, and we obtain closed
nilmanifolds carrying the same type of structure. Moreover one can see that
these MCP structures are not normal, their Reeb vector fields are however  Killing and hence $\h=0$.

\end{example}

\section{The tensor $\h$ and the Levi-Civita connection for MCP's}
Here we show some properties of the tensor field $\h$ for MCP manifolds. We also prove some formulae concerning the Levi-Civita connection $\nabla$ for a metric associated to a contact pair.

\subsection{The covariant derivative of $\phi$}\label{covariant-derivative-phi}

\begin{proposition}\label{pr:gen-metric}
Let $(\alpha_1, \alpha_2, \phi)$ be a contact pair structure
together  with a compatible metric $g$, and $\Phi$ the two-form defined by $\Phi (X,Y)= g(\phi X, Y)$. Then the covariant derivative of $\phi$ is given by
\begin{equation}\label{propert:gen-metric}
\begin{split}
2 g \left ((\nabla_X \phi)Y, W \right ) = 3 d\Phi (X,Y,W)-3 d\Phi(X, \phi Y, \phi W)+ g \left (N^1 (Y,W), \phi X \right )\\
+ 2 \sum_{i=1} ^2 \bigl (d\alpha_i (\phi Y , X) \alpha_i (W)- d\alpha_i (\phi W,X) \alpha_i (Y)\bigr )+ \sum_{i=1} ^2 \alpha_i (X) N^2_i (Y, W).
\end{split}
\end{equation}
 \end{proposition}

\begin{proof}
Applying the definition of the Levi-Civita connection to $2 g(\nabla _X  Y, W)$ and using the formula for the exterior derivative of $\Phi$ in terms of Lie brackets, we have:
\smallskip
\begin{equation*}
\begin{split}
&2 g((\nabla_X \phi)Y, W)\\
=\, \, &  2 g(\nabla_X (\phi Y) , W) +2 g(\nabla_X Y, \phi W)\\
=\, \,  & X \Phi (Y,W)+ \phi Y \left(\Phi (X, \phi W)+ \sum_{i=1} ^2 \alpha_i(X) \alpha_i (W)\right) + W \Phi (X,Y)\\
&+ \Phi ([X, \phi Y], \phi W)+\sum_{i=1} ^2 \alpha_i([X, \phi Y]) \alpha_i (W)- \Phi ([W,X], Y)-g (\phi [\phi Y , W], \phi X) \\
&- \sum_{i=1} ^2 \alpha_i([\phi Y, W]) \alpha_i (X)- X \Phi (\phi Y, \phi W) + Y \Phi (W,X) \\
&- \phi W \left(\Phi (X, \phi Y)+ \sum_{i=1} ^2 \alpha_i(X) \alpha_i (Y)\right)- \Phi ([X,Y], W) + \Phi ([\phi W, X], \phi Y)\\
&+ \sum_{i=1} ^2 \alpha_i([\phi W , X]) \alpha_i (Y)- g(\phi [Y, \phi W ], \phi X) + \sum_{i=1} ^2 \alpha_i(X) \alpha_i ( [\phi W , Y])\\
&\smallskip \\
&- g ([Y, W], \phi X) - \Phi ([Y,W], X)+ g([\phi Y , \phi W], \phi X) + \Phi ([\phi Y , \phi W], X)\\
&\smallskip \\
&+ g \left (2 d\alpha_1 (Y, W) Z_1 , \phi X \right ) + g (2 d\alpha_2 (Y, W) Z_2, \phi X)\\
&\smallskip \\
=\, \, &  3 d\Phi (X,Y,W)-3 d\Phi(X, \phi Y, \phi W)+ g \left (N^1 (Y,W), \phi X \right )\\
&+ 2 \sum_{i=1} ^2 \bigl (d\alpha_i (\phi Y , X) \alpha_i (W)- d\alpha_i (\phi W,X) \alpha_i (Y)\bigr )+ \sum_{i=1} ^2 \alpha_i (X) N^2_i (Y, W).
\end{split}
\end{equation*} \end{proof}
Applying Proposition \ref{pr:gen-metric} to a MCP with decomposable $\phi$, we obtain:

\begin{corollary}\label{cor1:gen-metric}
For an MCP $(\alpha_1, \alpha_2, \phi, g)$ with decomposable $\phi$, the covariant derivative of $\phi$ is given by
\begin{equation}\label{eq:cor1:gen-metric}
2 g((\nabla_X \phi)Y, W)= g \left (N^1 (Y,W), \phi X \right )+ 2 \sum_{i=1} ^2 \bigl (d\alpha_i (\phi Y , X) \alpha_i (W)- d\alpha_i (\phi W,X) \alpha_i (Y)\bigr ) .
\end{equation}
\end{corollary}

\begin{proof}
For an MCP with decomposable $\phi$ we have $N^2_1=N^2_2=0$ by Proposition \ref{p:Ni=0}. Moreover $-\Phi=d\alpha_1 + d\alpha_2$. Then \eqref{propert:gen-metric}
reduces to \eqref{eq:cor1:gen-metric}. \end{proof}
\begin{corollary}\label{cor2:gen-metric}
For an MCP with decomposable $\phi$ and Reeb vector fields $Z_1$~, $Z_2$
we have:
\begin{equation*}
\nabla _{Z_1} \phi =\nabla _{Z_2} \phi =0.
\end{equation*}
\end{corollary}
\begin{proof}
In \eqref{eq:cor1:gen-metric}, we put $X=Z_i$ for $i=1,2$, and we obtain $g((\nabla _{Z_i} \phi) Y , W)=0$.
\end{proof}

\subsection{The tensor field $\h$}\label{subsect-h}

When $\phi$ is decomposable so is the tensor field $\h$, because
for every $X\in \Gamma (T\mathcal{F}_i)$ we have $[Z_j , X] \in \Gamma (T\mathcal{F}_i)$ for $i,j=1,2$. In this case we have:
\begin{equation*}
\h=\h_1 \oplus \h_2 \; \text{and} \; \phi=\phi_1 \oplus \phi_2 \, ,
\end{equation*}
where $\h_1$ (respectively $\phi_1$) is the endomorphism of $ T\mathcal{F}_2$ induced by $\h$ (respectively by $\phi$) and vice-versa.
We can now state the following results:
\begin{theorem}\label{th:properties-h}
Let $(\alpha_1, \alpha_2, \phi, g)$ be an MCP with decomposable $\phi$ on a manifold $M$. Let
 $Z_1 , Z_2$ be the Reeb vector fields of $(\alpha_1, \alpha_2)$ and $Z=Z_1+Z_2$. Then we have:
\begin{enumerate}[{\indent(a)}]
\item
$\LL_{Z_1} \phi$~, $\LL_{Z_2} \phi$~, $\h$~, $\h_1$ and $\h_2$ are symmetric operators;
\item  $\nabla_X Z = -\phi X -\phi \h X $ for every $X \in \Gamma (TM)$;
\item $\h \circ \phi + \phi \circ \h= 0$ and $\h_i \circ \phi_i + \phi_i \circ \h_i= 0$ for $i=1,2$;
\item $\tr \h=\tr \h_1=\tr \h_2=0$~;
\item  $\alpha_i \circ \h = \alpha_i \circ \h _j =0$~ for every $i,j=1,2$.
\end{enumerate}
\end{theorem}
To prove this
we need the following:
\begin{lemma}\label{lemma-th:properties-h}
Let $(\alpha_1, \alpha_2, \phi, g)$ be an MCP  on a manifold $M$
with Reeb vector fields
 $Z_1$~and $Z_2$. Then for every $X \in \Gamma (TM)$, $\nabla _ X Z_1$ and $\nabla _ X Z_2$ are both tangent
 to the kernels of $\alpha_1$ and $\alpha_2$~.
\end{lemma}
\begin{proof}
Since $\nabla_{Z_i}Z_j=0$, it is enough to take $X$ horizontal.  Note also that
$\alpha_1(\nabla_XZ_2)=-\alpha_2(\nabla_XZ_1)$ and 
$\alpha_1(\nabla_XZ_1)=\alpha_2(\nabla_XZ_2)=0$. Now
$$\alpha_1(\nabla_XZ_2)=g(\nabla_{Z_2}X+[X,Z_2],Z_1)
=\alpha_1([X,Z_2])+Z_2\alpha_1(X)-g(X,\nabla_{Z_2}Z_1)=0$$
since $d\alpha_1(X,Z_2)=0$.
\end{proof}

\medskip
\begin{proof}[Proof of Theorem \ref{th:properties-h}]
~

\begin{enumerate}[{\noindent(a)}]
\item We want to show that $g \left (X, (\LL_{Z_j} \phi) Y \right )= g \left ((\LL_{Z_j} \phi) X, Y \right )$, for $j=1,2$. We prove the property for $j=1$, since the other case is similar.
For $X=Z_i$~, $i=1,2$ we have $g \left ((\LL_{Z_1} \phi) Z_i,  Y \right )= 0$ and $g \left (Z_i,(\LL_{Z_1} \phi) Y \right )=0$. The same holds for $Y=Z_i$~, $i=1,2$. Then we have to prove the symmetry of $\LL_{Z_1} \phi$ for $X, Y$ tangent to
$\ker \alpha_1 \cap \ker \alpha_2$.
By Corollary \ref{cor2:gen-metric} we have $\nabla_{Z_i} \phi=0$ for $i=1,2$. For $X, Y \in \ker \alpha_1 \cap \ker \alpha_2$, we have:
\begin{alignat*}{1}
g \left ((\LL_{Z_1} \phi )X, Y \right )= & g \left (- \nabla _{\phi X} Z_1+ \phi (\nabla_X Z_1), Y \right ) \\
= & g (Z_1 , \nabla_{\phi X} Y)- g (\nabla _X Z_1, \phi Y)\\
= & \alpha_1 (\nabla _{\phi X}Y)+ \alpha_1 (\nabla_X \phi Y) \\
= & \alpha_1 (\nabla _ Y \phi X)+ \alpha_1 (\nabla_{\phi Y} X) \\
= & g \left (X, (\LL_{Z_1} \phi ) Y \right ) ,
\end{alignat*}
where we have used that $Z_1$ is orthogonal to $X, Y \in \ker \alpha_1 \cap \ker \alpha_2$ and that for an MCP the tensors $N^2_1$ and $N^2_2$ vanish.
It is clear that $\h$ is symmetric as well and after restriction this is also true for $\h_1$ and $\h_2$.

\item By Corollary \ref{cor1:gen-metric}, for $i=1,2$, and for every $X, Y \in \Gamma (TM)$, we have:
\begin{alignat*}{1}
2 g \left ((\nabla_X \phi)Z_i, Y \right )= & g \left (N^1 (Z_i , Y), \phi X \right )- 2 d\alpha_i (\phi Y,X) \\
= & g \left ( \phi^2 [Z_i , Y] - \phi [Z_i , \phi Y], \phi X \right )- 2 d\alpha_i (\phi Y,X) \\
= & - g \left ( \phi (\LL_{Z_i} \phi)Y), \phi X \right )- 2 d\alpha_i (\phi Y,X) \\
= & -  g \left ( (\LL_{Z_i} \phi)Y), X \right )+ \Big(\sum_{j=1} ^2 \alpha_j( \left (\LL_{Z_i} \phi)Y \right )
\alpha_j (X) \Big)- 2 d\alpha_i (\phi Y,X) \\
= & -  g \left ( (\LL_{Z_i} \phi)Y, X \right )- 2 d\alpha_i (\phi Y,X) .
\end{alignat*}
Then we obtain:
\begin{equation*}
\begin{split}
\sum_{i=1} ^2  2g \left ((\nabla_X \phi)Z_i, Y \right )= & \sum_{i=1} ^2 \bigl ( -2 g \left( (\LL_{Z_i} \phi) Y ,  X \right )- 2 d\alpha_i (\phi Y,X) \bigr )\\
=& -2 \sum_{i=1} ^2 g \left ( (\LL_{Z_i} \phi) Y ,  X \right ) - 2 g (\phi Y, \phi X) \\
=& -2 \sum_{i=1} ^2 g \left ( (\LL_{Z_i} \phi )Y ,  X \right ) - 2 g (Y,  X) + 2 \sum_{i=1} ^2 \alpha_i (X) \alpha_i (Y)\\
\end{split}
\end{equation*}
and then
\begin{equation*}
\begin{split}
g((\nabla_X \phi)Z, Y)= &  - g ( \h Y , X)- g (X , Y) + \sum_{i=1} ^2 \alpha_i (X) \alpha_i (Y) \\
=&  - g ( \h Y , X)- g (X , Y) + \sum_{i=1} ^2 g \left (\alpha_i (X) Z_i ,Y \right ).
\end{split}
\end{equation*}
Since the last equation is true for every $X, Y \in \Gamma (TM)$, this implies
\begin{equation*}
- \phi \nabla _ X Z = (\nabla_X \phi)Z= - \h X - X  + \alpha_1 (X) Z_1 + \alpha_2 (X) Z_2.
\end{equation*}
Applying $\phi$ to the last equation and using Lemma \ref{lemma-th:properties-h} gives:
\begin{equation*}
\nabla_X Z=  - \phi X  - \phi \h X+ \alpha_1 (\nabla_X Z) Z_1 + \alpha_2 (\nabla_X Z) Z_2= - \phi X  - \phi \h X.
\end{equation*}

\item For $X, Y \in \Gamma(TM)$, by the symmetry  of $\h$ and the formula $\nabla_X Z=- \phi X  - \phi \h X$, we have:
\begin{alignat*}{1}
2 g (X, \phi Y)=& 2 (d\alpha_1 + d\alpha_2) (X, Y) \\
= & \sum_{i=1} ^2 \bigl (g(\nabla_X Z_i , Y) -g(\nabla_Y Z_i , X) \bigr ) \\
= & g(\nabla_X Z , Y) -g(\nabla_Y Z , X) \\
= & - g (\phi X , Y)+ g (\phi Y , X)- g (\phi \h X , Y)+ g (\phi \h Y , X) \\
= & - g (\phi \h X , Y)+ g (\phi \h Y , X) + 2 g (X, \phi Y )\\
=& g (\h \phi Y+\phi \h Y , X)+ 2 g (X, \phi Y ),
\end{alignat*}
which implies that $\h \circ \phi +\phi \circ  \h=0$. After restriction of $\h$ and $\phi$ to $T\mathcal{F}_i$ for $i=1,2$, we obtain $\h_i \circ  \phi_i +\phi_i \circ  \h_i=0$.

\item Since the endomorphism $\h$ is symmetric, at every point $p \in M$ there exists an eigenbasis of $T_p M$. Let $V$ be
 an eigenvector relative to the eigenvalue $\lambda$. Then, by (c), we have:
$$
\h_p (\phi _p V)= - \lambda (\phi _p V),
$$
which means that $-\lambda$ is also an eigenvalue, relative to the eigenvector  $\phi _p V$, and then the trace of $\h_p$ vanishes for every $p \in M$. Similarly we have $\tr \h_1=\tr \h_2=0$.

\item The last property follows easily from (c).
\end{enumerate}
\end{proof}

\begin{corollary}
Let $(\alpha_1, \alpha_2, \phi, g)$ be an MCP with decomposable $\phi$ and Reeb vector fields $Z_1$~, $Z_2$~.
If the vector field $Z=Z_1 + Z_2$ is Killing,
 then we have
\begin{equation*}
\nabla_X Z = -\phi X.
\end{equation*}
\end{corollary}
\begin{proof}
By Proposition \ref{p:Z-killing}, the vector field $Z$ is Killing if and only if $\h=0$. Applying this to Theorem \ref{th:properties-h}-(b), we get $\nabla_X Z = -\phi X$.
\end{proof}
A special case is given when both Reeb vector fields are Killing. A first example of the latter situation concerns the non-normal MCP on the nilpotent Lie group $G$ and its closed nilmanifolds $G/\Gamma$ described in Example \ref{ex:nilp-lie-group}.

One can also have $Z_i$ Killing by choosing normal structures (see Corollary \ref{cor:normal-killing}). Then
here is a second example, with a normal MCP but where the manifold is not a product of two metric contact manifolds:

\begin{example}
Let $M=\widetilde{SL_2}$ be the universal covering of the identity
component of the isometry group of the hyperbolic plane $\bH^2$
endowed with an invariant Sasakian structure $(\alpha, \phi, g)$ (see
\cite{Geiges2}) and $N= M \times M$. It is well known that $N$ admits cocompact irreducible lattices $\Gamma$ (see
\cite{Bor}). This means that $ \Gamma$ does not admit any subgroup
of finite index which is a product of two lattices of $M$. The
manifold $N$ can be endowed with the product MCP
structure and by the invariance of the structure by $\Gamma$,
the MCP descends to the quotient and is normal. Even if
the local structure is like a product, globally the two characteristic foliations can
be very interesting in the sense that both could have dense
leaves.
\end{example}

\section{Some curvature properties}
In this section, for a manifold $M$ carrying an MCP $(\alpha_1, \alpha_2, \phi, g)$ with decomposable $\phi$, we set $Z=Z_1 + Z_2$ where $Z_1$, $Z_2$ are the Reeb vector fields.
We prove some properties of the curvature tensor and the Ricci curvature, which are analogous to those of metric contact structures (see e.g. \cite{Blairbook}).
As a consequence we prove the non-flatness of metrics associated to normal MCP's. This implies the non-existence of flat non-K\"ahler Vaisman manifolds.

\subsection{The curvature}\label{s:curvature}
We denote by
$\Curv$
the curvature tensor  of the metric $g$, and by $\Ric$ its Ricci curvature.

\begin{proposition}
Let $(\alpha_1, \alpha_2, \phi, g)$ be an MCP with decomposable $\phi$ on a manifold $M$. Then:
\begin{alignat}{1}
&(\nabla_Z \h) X= \phi X- \h^2\phi X - \phi \Curv_{X Z} Z \label{th-curv-prop:eq1}\\
&\frac{1}{2} \bigl (\Curv_{Z X}Z - \phi (\Curv_{Z \phi X}Z) \bigr )= \phi^2 X + \h^2 X \label{th-curv-prop:eq2}.
\end{alignat}
\end{proposition}

\begin{proof}
Corollary \ref{cor2:gen-metric} implies $\nabla_Z (\phi X)= \phi (\nabla _Z X)$. Using this and Theorems \ref{th:ass-metrics} and \ref{th:properties-h},  we apply $\phi$ to:
\begin{equation*}
R _{Z X} Z= \nabla _Z \nabla_ X Z -\nabla _X \nabla_ Z Z - \nabla _{[Z , X]} Z = \nabla _ Z (- \phi X - \phi \h X) + \phi [Z, X]+ \phi \h [Z,X] ,
\end{equation*}
and we obtain:
\begin{alignat*}{1}
\phi(R _{ZX}Z)=&\nabla _ Z ( X + \h X) - \alpha_1 \bigl (\nabla _ Z ( X + \h X) \bigr ) Z_1- \alpha_2 \bigl (\nabla _ Z ( X + \h X) \bigr ) Z_2- [Z,X]\\
+&\alpha_1 ([Z,X]) Z_1+\alpha_2 ([Z,X]) Z_2 - \h [Z,X]+ \alpha_1 (\h[Z,X]) Z_1+\alpha_2 (\h[Z,X]) Z_2\\
=& (\nabla_Z \h)X + \nabla _X Z + \h \nabla _ X Z \\
= & (\nabla_Z \h)X - \phi X - \phi \h X + \h (-\phi X - \phi \h X)\\
=& (\nabla_Z \h)X -\phi X + \h^2 \phi X ,
\end{alignat*}
which gives \eqref{th-curv-prop:eq1}.

To prove \eqref{th-curv-prop:eq2} we first remark that $R_{Z X} Z$ is tangent to the kernels of $\alpha_1$ and $\alpha_2$. Then we have:
$$
R_{Z X} Z=-\phi^2 R_{Z X} Z= \phi^2 X - \phi \h^2 \phi X - \phi \bigl ((\nabla_Z \h) X \bigr )= \phi^2 X+ \h ^2 X- \phi \bigl ((\nabla_Z \h) X \bigr ) .
$$
Using the previous expression and taking the difference $R_{Z X} Z - \phi (R_{Z \phi X} Z)$ gives \eqref{th-curv-prop:eq2}.
\end{proof}

\begin{theorem}\label{th:ricci-curvature1}
Let $(\alpha_1, \alpha_2, \phi, g)$ be an MCP of type $(h,k)$ with decomposable $\phi$ on a $(2h+2k+2)$-dimensional manifold $M$. Then we have:
\begin{equation}\label{th:ricci-curvature1-eq1}
\Ric (Z)= h+k - \frac{1}{2}\tr \h^2 .
\end{equation}
Moreover $\Ric (Z)= h+k $ if and only if $Z$ is Killing.
\end{theorem}
\begin{proof}
Denote by $K (Z, X)$ the sectional curvature of the plane determined by $\{Z, X\}$. By using \eqref{th-curv-prop:eq2} for $X$ of unit length and orthogonal to $Z_1$ and $Z_2$, and recalling that $g(Z, Z)= 2 $, we obtain
\begin{alignat*}{1}
\Sec (Z, X)+ \Sec (Z, \phi X)= & - \frac{1}{2}g \bigl ( \Curv_{Z X}Z - \phi (\Curv_{Z \phi X}Z), X \bigr ) \\
= &- g( \phi^2 X + \h^2 X , X) \\
= & 1 - g (\h^2 X , X) . \\
\end{alignat*}
Let $\{ Z_1 , Z_2, X_1, \cdots ,X _{2h+2k} \}$ be a local $\phi$-basis, that is an orthogonal basis for which the $X_i$ have unit length and $X_{2i}=\phi X_{2i-1}$. Then, since $K \left (Z,Z_1-Z_2 \right )=0$,
taking the sum  $\sum_{i=1}^{2h+2k} \Sec (Z, X_i)$ we obtain \eqref{th:ricci-curvature1-eq1}. Now $\Ric (Z)= h+k $ if and only if $\tr \h^2=0$. Because $\h$ is symmetric the trace of $\h^2$ vanishes if and only if $\h=0$.
Now use Proposition \ref{p:Z-killing} to complete the proof.
\end{proof}

The following result generalizes to MCP's a theorem of Hatakeyama et al. \cite{HOT}:

\begin{theorem}\label{th:tanno}
Let $(\alpha_1, \alpha_2, \phi, g)$ be an MCP with decomposable $\phi$ on a $(2h+2k+2)$-dimensional manifold $M$. Then $Z$ is Killing if and only if for all the plane sections $(Z, X)$ with $X$ orthogonal to both $Z_1$ and $Z_2$, the value of the sectional curvature $\Sec (X, Z)$ is $1/2$. Moreover in this case, for every $Y$ we have:
\begin{equation}\label{th:tanno-eq}
\Curv _{Y Z} Z= Y - \alpha_1 (Y) Z_1 - \alpha_2 (Y) Z_2~.
\end{equation}
\end{theorem}

\begin{proof}
If for all the plane sections $(Z, X)$ with $X$ orthogonal to $Z_1$ and  $Z_2$, we have $\Sec (X, Z)=1/2$. Then $\Ric(Z)= h+k$,
and $\h =0$ by Theorem \ref{th:ricci-curvature1}~.

Conversely, let $\h =0$. Using $\nabla_X Z= -\phi X$, for $X$ of unit length and orthogonal to $Z_1$ and  $Z_2$, and recalling
that $\nabla_Z Z=0$, we have
\begin{alignat*}{1}
2 \Sec (Z, X)= & - g ( R_{Z X} Z ,X)\\
= & - g (\nabla _ Z \nabla_ X Z - \nabla _{[Z, X] } Z , X)\\
= & g ( \nabla _ Z \phi X - \phi [Z, X], X)\\
= & g (\phi ( \nabla _ Z X )- \phi [Z, X], X)\\
= & g (\phi ( \nabla _ X Z + [Z, X] )- \phi [Z, X], X)\\
= & g (\phi ( \nabla _ X Z  ), X)\\
= & g (-\phi^2 X , X)\\
= & 1.
\end{alignat*}

To obtain \eqref{th:tanno-eq}, we have just to set $\h =0$ in \eqref{th-curv-prop:eq1}, apply $\phi$ on both sides and observe that an easy computation gives $g(R_{Y Z} Z,Z_i)=0$.
\end{proof}

\subsection{Normal MCP's and Vaisman structures}\label{s:vaisman}

By Proposition \ref{prop:normal-lziphi-vanish}, for a normal contact pair the tensor $\h$ vanishes necessarily. Thus, by \eqref{th:ricci-curvature1-eq1}, we have:

\begin{corollary}
A metric associated to a normal contact pair with decomposable endomorphism cannot be flat.
\end{corollary}

In particular this is true for normal MCP's of type $(h,0)$ which are nothing but non-K\"ahler Vaisman structures modulo constant rescaling of the metric \cite{BK2}. For this case the previous Corollary can be stated as:

\begin{theorem}
The metric of a non-K\"ahler Vaisman
manifold cannot be flat.
\end{theorem}

Here compactness is not needed. However this result was known for closed manifolds. Indeed a Vaisman structure is a particular locally conformally K\"ahler (lcK) structure. According to \cite{Vaisman} (see also \cite[Proposition 2.5]{Ornea}) a closed lcK manifold of constant curvature is necessary K\"ahler. Hence flat non-K\"ahler Vaisman structures do not exist on closed manifolds.

\subsection{}

In complete analogy to the case of contact metric manifolds, we want to define two tensor fields that are useful for the calculations in the problem of finding metric contact pairs with curvature vanishing on the vertical subbundle. First observe that for a metric contact pair with decomposable $\phi$, we have:
\begin{equation}\label{eq:*}
\begin{split}
&2g((\nabla_X \phi)W, \phi Y )-2g((\nabla_X \phi)\phi W, Y )\\
=&g(N^1(W ,\phi Y)-N^1(\phi W , Y), \phi X)-2d\alpha_1(\phi^2 Y, X)\alpha_1(W)-2d\alpha_2(\phi^2 Y, X)\alpha_2(W)\\
&-2d\alpha_1(\phi^2 W, X)\alpha_1(Y)-2d\alpha_2(\phi^2 W, X)\alpha_2(Y)\\
=&\alpha_1(Y) g([\phi W, Z_1]-\phi [W, Z_1], \phi X) +\alpha_2(Y)g([\phi W, Z_2]-\phi [W, Z_2], \phi X)\\
&+\alpha_1(W) g([\phi Y, Z_1]-\phi [Y, Z_1], \phi X) +\alpha_2(W)g([\phi Y, Z_2]-\phi [Y, Z_2], \phi X)\\
&+2d\alpha_1(Y,X)\alpha_1(W)+2d\alpha_2(Y,X)\alpha_2(W)+2d\alpha_1(W,X)\alpha_1(Y)+2d\alpha_2(W,X)\alpha_2(Y)\, .
\end{split}
\end{equation}

\bigskip

Replacing $W$ with $\phi W$ in \eqref{eq:*}, we obtain the following:
\begin{equation}\label{eq:**}
\begin{split}
&2g((\nabla_X \phi)\phi W, \phi Y )+2g((\nabla_X \phi)W, Y )-2g(Y, (\nabla_X\phi)(\alpha_1(W)Z_1 + \alpha_2(W)Z_2))\\
=&\alpha_1(Y) g(-[ W, Z_1]-\phi [\phi W, Z_1], \phi X) +\alpha_2(Y)g(-[W, Z_2]-\phi [\phi W, Z_2], \phi X)\\
&+2d\alpha_1(\phi W,X)\alpha_1(Y)+2d\alpha_2(\phi W,X)\alpha_2(Y) \, .
\end{split}
\end{equation}

Taking $X,Y,W$ horizontal in \eqref{eq:*} and  in \eqref{eq:**}, we have:

\begin{equation}\label{eq:*hor}
g((\nabla_X \phi)W, \phi Y )=g((\nabla_X \phi)\phi W, Y ) \, .
\end{equation}

\begin{equation}\label{eq:**hor}
g((\nabla_X \phi)\phi W, \phi Y )=-g((\nabla_X \phi)W, Y ) \, .
\end{equation}

Using (4) with $X, Y,W$ horizontal we get
\begin{equation}\label{eq:***}
g((\nabla _X\phi)Y, W)+g((\nabla _{\phi X}\phi)\phi Y, W)=0 \, .
\end{equation}

For the curvature operator we have:
\begin{equation}
\begin{split}
\Curv _{X Y}Z&=-\nabla_X(\phi Y+\phi \h Y)+\nabla_Y(\phi X+\phi \h X)+\phi [X,Y]+\phi \h [X,Y]\\
&=-(\nabla_X \phi)Y+(\nabla_Y \phi)X-(\nabla_X \phi\h)Y+(\nabla_X \phi\h)Y \, ,
\end{split}
\end{equation}
which gives:
\begin{equation}\label{eq:curv1}
g(\Curv _{Z W}X,Y)=-g(X, (\nabla_W \phi)Y)-g(W,(\nabla_X \phi\h)Y)+g(W,(\nabla_Y \phi\h)X) \, ,
\end{equation}
or equivalently
\begin{equation}
g(\Curv _{Z X}Y,W)=-g(Y, (\nabla_X \phi)W)-g(X,(\nabla_Y \phi\h)W)+g(X,(\nabla_W \phi\h)Y) \, .
\end{equation}
Now we define the following tensors:
\begin{definition}
For $X,Y,W \in \Gamma(TM)$, set
\begin{equation}\label{eq:def-A}
\begin{split}
A(X,Y,W)=&-g(Y, (\nabla_X \phi)W)+g(\phi Y, (\nabla_X \phi)\phi W)\\
&-g(Y, (\nabla_{\phi X} \phi)\phi W)-g(\phi Y, (\nabla_{\phi X} \phi)W) \, ,
\end{split}
\end{equation}
\begin{equation}\label{eq:def-B}
\begin{split}
B(X,Y,W)=&-g(X, (\nabla_Y \phi\h)W)+g(X, (\nabla_{\phi Y} \phi\h)\phi W)\\
&-g(\phi X, (\nabla_Y \phi\h)\phi W)-g(\phi X, (\nabla_{\phi Y} \phi\h)W) \, .
\end{split}
\end{equation}
\end{definition}
Taking $X,Y,W$ horizontal and using (16), we obtain the following relation:
\begin{equation}\label{eq:rel-AB}
\begin{split}
A(X,Y,W)+B(X,Y,W)-B(X,W,Y)=&g(\Curv _{Z X}Y,W)-g(\Curv _{Z X}\phi Y,\phi W)\\
& +g(\Curv _{Z \phi X}Y,\phi W)+g(\Curv _{Z \phi X}\phi Y,W) \, .
\end{split}
\end{equation}

The following lemma will be useful in the sequel:
\begin{lemma}\label{l:propA-B-horizontal}
For every $X,Y,W$ horizontal, we have:
$$
A(X,Y,W)+B(X,Y,W)-B(X,W,Y)=-2g((\nabla_{\h X} \phi)Y ,W)
$$
\end{lemma}
\begin{proof}
For $X,Y,W$ horizontal, using \eqref{eq:*hor},  \eqref{eq:**hor} and \eqref{eq:***} we have:
\begin{equation}
\begin{split}
A(X,Y,W)=&-2g(Y, (\nabla_X \phi)W)-2g(Y, (\nabla_{\phi X} \phi)\phi W)=0\, .
\end{split}
\end{equation}
Also for $X,Y,W$ horizontal, we calculate
\begin{equation}
\begin{split}
B(X,Y,W)=&-g(X, (\nabla_Y \phi)\h W +\phi (\nabla_Y \h)W)+g(X, (\nabla_{\phi Y} \h) W -\phi \h(\nabla_{\phi Y} \phi)W)\\
&-g(\phi X, (\nabla_Y \h) W -\phi \h(\nabla_Y \phi)W)-g(\phi X, (\nabla_{\phi Y} )\h W +\phi (\nabla_{\phi Y} \h)W)\\
=&-g(X, (\nabla_Y \phi)\h W )+g(X,  \h\phi(\nabla_{\phi Y} \phi)W)\\
&+g(\phi X, \phi \h(\nabla_Y \phi)W)-g(\phi X, (\nabla_{\phi Y}\phi )\h W)\, .
\end{split}
\end{equation}
Now, we  have
\begin{equation}
\begin{split}
-g(\phi \h X, (\nabla_{\phi Y} \phi)W)=&g((\nabla_{\phi Y} \phi)\phi \h X, W)\\
=&-g((\nabla_{Y} \phi)\h X, W)\\
=&g(\h X, (\nabla_{Y} \phi)W) \, ,
\end{split}
\end{equation}
where in the second line we used \eqref{eq:***}, and furthermore we also have
\begin{equation}
\begin{split}
-g(\phi X, (\nabla_{\phi Y} \phi)\h W)=&g((\nabla_{\phi Y}\phi )\phi X,\h W)\\
=&-g((\nabla_{Y} \phi) X,\h W)\\
=&g(X,(\nabla_{Y}\phi )\h W) \, ,
\end{split}
\end{equation}
again using \eqref{eq:***}.
This in turn gives
$$
B(X,Y,W)=2g(\h X,(\nabla_{Y} \phi) W) \, ,
$$
and putting all this together, we obtain
\begin{equation}
\begin{split}
A(X,Y,W)+B(X,Y,W)-B(X,W,Y)=&2g(\h X,(\nabla_{Y}\phi ) W)-2g(\h X,(\nabla_W \phi) Y)\\
=&-2g((\nabla_{\h X} \phi)Y ,W).
\end{split}
\end{equation}
\end{proof}

\begin{corollary}\label{c:rel-AB}
If the curvature of an MCP with decomposable $\phi$ satisfies $\Curv _{XY}Z_1=\Curv_{XY}Z_2=0$ for all $X,Y$, then for all horizontal vector fields
$$
g((\nabla_{\h X} \phi)Y ,W)=0.
$$
\end{corollary}
\begin{proof}
The left hand side of the equation of Lemma \ref{l:propA-B-horizontal} vanishes by the assumption on the curvature and by \eqref{eq:rel-AB}.
\end{proof}

\section{Curvature vanishing on the vertical subbundle}
In this section we prove for MCP's the analogues of the results of the second named author on metric contact manifolds \cite{Blair3, Blair4}. Recall that for a contact pair we defined the vertical subbundle as the subbundle $\mathcal{V}$ spanned by the Reeb vector fields $Z_1, Z_2$. We will say that an MCP  has {\it curvature vanishing on the vertical subbundle} if the following condition is satisfied for all vector fields $X,Y$:
$$
\Curv _{XY} Z_i=0\, , \; \text{for} \, i=1,2\, .
$$ 
The standard example of such a situation is the product of the unit tangent bundles of $\bE^h$ and $\bE^k$, since each of them is endowed with a metric contact structure with curvature vanishing along the corresponding Reeb vector field (see \cite{Blair4}). Actually Theorem \ref{th:flat-metrics} below says exactly that locally this is the only possibility. 

Before the statement of the theorem, we make some remarks. First, observe that a contact pair of type $(0,0)$ is in some sense trivial, because we would like to have an induced contact form on the leaves of at least  one of the characteristic foliations. Moreover, if the manifold is endowed with an associated metric then it is flat and then locally isometric to $\bE^2$. Since the forms composing the contact pair play a symmetric role, to exclude the former trivial case we only consider contact pairs of type $(h,k)$ with $h\geq1$.

\begin{theorem}\label{th:flat-metrics}
Let $M$ be a $(2h+2k+2)$-dimensional manifold endowed with a metric contact pair $(\alpha_1, \alpha_2, \phi, g)$ of type $(h,k)$ (with $h\geq 1$) and decomposable $\phi$. If the curvature of $g$ vanishes on the vertical subbundle,
then $M$ is locally isometric to $\bE^{h+1}\times \bS^h(4)\times \bE^{k+1}\times \bS^k(4)$ if $k\geq 1$ or  to $\bE^{h+1}\times \bS^h(4)\times \bE^{1}$ if $k=0$.
\end{theorem}
\begin{proof}
We split the proof
into several steps:

\begin{noindenum}
\item[a)] We have seen that the decomposability of $\phi$ implies that $\h$ is also decomposable and we set $\phi=\phi_1 \oplus \phi_2 $ and $\h=\h_1 \oplus \h_2 $ as in Section \ref{subsect-h}. If the curvature tensor $\Curv$ vanishes on the vertical subbundle, by \eqref{th-curv-prop:eq2} we have $\h^2=- \phi ^2$ and then, since $\h$ is symmetric, for its rank  we have $\rk \h=\rk \h^2=\rk \phi^2=\rk \phi=2h+2k$. If $X$ is an eigenvector of $\h$ corresponding to a non-zero eigenvalue $\lambda$ then it is orthogonal to the Reeb vector fields (which are $0$-eigenvectors) and we have:
\begin{equation*}
\lambda ^2 g(X,X)=g(\lambda X, \lambda X)=g(\h X, \h X)=g(\h^2 X, X)=-g(\phi^2 X,X)=g(X,X) .
\end{equation*}
Thus the non-zero eigenvalues of $\h$ are $\pm 1$.  By restriction the same is true for the eigenvalues of $\h_1$ and $\h_2$. Moreover the $(\pm 1)$-eigenspaces of $\h$ are direct sums of the $(\pm 1)$-eigenspaces of $\h_1$ and $\h_2$ at every point. Also observe that $\phi$ (respectively $\phi_j$) intertwines the eigenspaces corresponding to $+1$ and $-1$, because it anticommutes with $\h$ (respectively $\h_j$).

Let $[-1]_1$, $[-1]_2$, $[-1]$ (respectively $[+1]_1$, $[+1]_2$, $[+1]$) be the $(-1)$-eigendistributions (respectively $(+1)$-eigendistributions) of $\h _1$, $\h_2$ and $\h$ respectively. By the previous discussion we have $[-1]=[-1]_1 \oplus [-1]_2$ and $[+1]=[+1]_1 \oplus [+1]_2$.

By Theorem \ref{th:properties-h}, we obtain
\begin{equation}\label{eq1-th:flat-metrics}
\nabla_X Z=0 \; ,\; \forall X \in [-1].
\end{equation}
Moreover, the vanishing of $\Curv$ on the vertical subbundle implies that for every vector fields $X,Y\in [-1]$ (respectively $[-1]_j$) we have:
\begin{equation*}
0=\Curv _{X Y} Z=-\nabla _{[X,Y]} Z =\phi [X,Y]+\phi \h [X,Y] \, ,
\end{equation*}
which implies $ \h \phi [X,Y]= \phi[X,Y]$ and then $\phi [X,Y] \in [+1]$ (respectively $[+1]_j$). Applying $\phi$ to both sides of the equation gives $\phi^2 [X,Y] \in [-1]$ (respectively $[-1]_j$).  Calculating further we have
$$
\phi^2 [X,Y]=-[X,Y]+d\alpha_1 (X,Y)Z_1 +d\alpha_2 (X,Y) Z_2= -[X,Y] \, .
$$
The last equation is clear if $X,Y$ are tangent to different foliations, since in this case the $d\alpha_i$ vanish, and is easily deduced from the following observation when $X,Y$ are tangent to the same foliation. In fact if $X,Y$ are tangent to the same foliation, say $\FF_1$ for example, we have $d\alpha_2 (X,Y)=d\alpha_1 (X,Y)+d\alpha_2 (X,Y)=g(X,\phi Y)=0 $ since $X\in [-1]$ and $\phi$ intertwines the $(\pm 1)$-eigenspaces.

The same calculations give $[X, Z_i]\in [-1]$ (respectively $[-1]_j$) for every $X\in [-1]$ (respectively $[-1]_j$) and $[\phi Y, Z_i]\in [-1]$ (respectively $[-1]_j$) for every $Y\in [+1]$ (respectively $[+1]_j$). In particular this implies that the distributions $[-1]_j $, $[-1]$, $[-1]_j \oplus \bR Z_i$, $[-1]\oplus \bR Z_i$, $[-1]_j \oplus \bR Z_1 \oplus \bR Z_2$  and $[-1]\oplus \bR Z_1 \oplus \bR Z_2$ are integrable.
\medskip
\item[b)]
According to the local model for contact pairs of type $(h,k)$, on every point there exist local coordinates
$
(u_0, \cdots , u_{2k}, v_0, \cdots , v_{2h})
$
such that
$\frac{\partial}{\partial u_0}, \cdots, \frac{\partial}{\partial u_{2k}}$
span
$\Ker \alpha_1 \cap \Ker d\alpha_1$
and
$\frac{\partial}{\partial v_0}, \cdots, \frac{\partial}{\partial v_{2h}}$
span
$\Ker \alpha_2 \cap \Ker d\alpha_2$.
By the integrability of $[-1]_2 \oplus \bR Z_2$ (respectively $[-1]_1 \oplus \bR Z_1$), these local coordinates can be chosen
such that
one also has that
$\frac{\partial}{\partial u_0}, \cdots, \frac{\partial}{\partial u_k}$ span $[-1]_2\oplus \bR Z_2$ and $\frac{\partial}{\partial v_0}, \cdots, \frac{\partial}{\partial v_h}$ span $[-1]_1 \oplus \bR Z_1$.

Let us define the vector fields $X_i=\frac{\partial}{\partial u_{k+i}}+\sum_{p=0}^k f_{ip}\frac{\partial}{\partial u_p}$, for $1\leq i\leq k$ and $Y_s=\frac{\partial}{\partial v_{h+s}}+\sum_{q=0}^h \tilde f_{sq}\frac{\partial}{\partial v_q}$ for $1\leq s\leq h$, where the functions $f_{ip}$~,  $\tilde f_{sq}$ are chosen in such a way that $X_i\in [+1]_2$ and $Y_s \in [+1]_1$. In general those functions depend on all coordinates. It is clear that at every point of $M$, the $X_i$'s and $Y_s$'s form a basis for $[+1]_2$ and $[+1]_1$ respectively.

A direct calculation with local coordinates gives:
\begin{equation*}
\begin{split}
&\Big[\frac{\partial}{\partial u_p}, X_i \Big] \in [-1]_2 \oplus \bR Z_2 \; ,\;  0\leq p \leq k \; ,\;  0\leq i \leq k  \; , \\
&\Big[\frac{\partial}{\partial v_q}, Y_s\Big]  \in [-1]_1 \oplus \bR Z_1  \; , \;0 \leq q \leq h \; ,\;  0\leq s \leq h \; , \\
&\Big[\frac{\partial}{\partial u_p}, Y_s\Big]  \in [-1]_1 \oplus \bR Z_1 \; ,\;  0\leq p \leq k \; ,\;  0\leq s \leq h \; , \\
&\Big[\frac{\partial}{\partial v_q}, X_i\Big] \in [-1]_2 \oplus \bR Z_2  \; , \;0 \leq q \leq h \; ,\;  0\leq i \leq k \; , \\
\end{split}
\end{equation*}
and
\begin{equation*}
[X_i, Y_s] \in [-1] \oplus \bR Z_1 \oplus \bR Z_2\,.
\end{equation*}
Then we have
\begin{equation*}
\nabla_{[\frac{\partial}{\partial u_p}, X_i ] }Z=\nabla_{[\frac{\partial}{\partial u_p}, Y_s ] }Z=\nabla_{[\frac{\partial}{\partial v_q}, X_i ] }Z=\nabla_{[\frac{\partial}{\partial v_q}, Y_s ] }Z =\nabla_{[X_i, Y_s] }Z =0\, .
\end{equation*}
The assumption on the curvature implies $\Curv _{\frac{\partial}{\partial u_p} X_i } Z=0$,
$0\leq p\leq k$, and we obtain
\begin{equation}\label{eq1.5-th:flat-metrics}
0=\nabla_{\frac{\partial}{\partial u_p} }\nabla _{X_i} Z - \nabla _{X_i}  \nabla_{\frac{\partial}{\partial u_p} }Z=-2\nabla_{\frac{\partial}{\partial u_p} } \phi X_i \, .
\end{equation}
Since the $\frac{\partial}{\partial u_p}$'s span $[-1]_2\oplus \bR Z_2$ and the connection is tensorial in the first argument, we have:
\begin{equation}\label{eq2-th:flat-metrics}
\nabla_{\phi X_j} {\phi X_i}=0 \; \; , \forall i,j \, .
\end{equation}
In a similar way we obtain the following formulae:
\begin{equation}\label{eq3-th:flat-metrics}
\begin{split}
\nabla_{\phi Y_r} {\phi Y_s}=0 & \; \; , \forall r,s \; ,\\
\nabla_{\phi X_i} {\phi Y_r}=0 & \; \; , \forall i,r \; ,\\
\nabla_{\phi Y_r} {\phi X_i}=0 & \; \; , \forall i,r \, .
\end{split}
\end{equation}
These imply that the integral submanifolds of $ [-1] \oplus \bR Z_1 \oplus \bR Z_2$ are totally geodesic and flat.

A direct calculation with local coordinates shows that $[X_i, X_j]$ is in $[-1]_2\oplus \bR Z_2$. Differentiating $g(X_i , Z)=0$ along $X_j$  we obtain $g(\nabla_{X_j}X_i , Z)=0$. Interchanging $i$ and $j$ and taking the difference we get $0=g([X_i , X_j], Z)=g([X_i , X_j], Z_2)$ since   $[X_i, X_j]$ is orthogonal to $Z_1$. This actually means that $[X_i, X_j]$ is in $[-1]_2$.

Then we have
\begin{equation*}
0= \Curv _{X_i X_j}Z=-2 \nabla _{X_i} \phi X_j+ 2 \nabla _{X_j} \phi X_i \, ,
\end{equation*}
or equivalently
\begin{equation}\label{eq3.5-th:flat-metrics}
 \nabla _{X_i} \phi X_j=\nabla _{X_j} \phi X_i \, .
\end{equation}
Similarly we obtain
\begin{equation}\label{eq4-th:flat-metrics}
 \nabla _{Y_r} \phi Y_s=\nabla _{Y_s} \phi Y_r \, .
\end{equation}
With similar calculations, using the fact that $[X_i , Y_r] \in [-1] \oplus \bR Z_1 \oplus \bR Z_2$, we obtain
\begin{equation}\label{eq5-th:flat-metrics}
 \nabla _{X_i} \phi Y_r=\nabla _{Y_r} \phi X_i \, .
\end{equation}
Equations \eqref{eq3.5-th:flat-metrics}--\eqref{eq5-th:flat-metrics} can also be written as follows:
\begin{equation}\label{eq6-th:flat-metrics}
\begin{split}
\phi[X_i, X_j]&=-( \nabla _{X_i} \phi) X_j+(\nabla _{X_j} \phi )X_i \, ,\\
\phi[Y_r, Y_s]&=-( \nabla _{Y_r} \phi) Y_s+(\nabla _{Y_s} \phi )Y_r \, ,\\
\phi[X_i, Y_r]&=-( \nabla _{X_i} \phi) Y_r+(\nabla _{Y_r} \phi )X_i \, .\\
\end{split}
\end{equation}
Using \eqref{eq1-th:flat-metrics} and \eqref{eq2-th:flat-metrics} we obtain
$$
0=\Curv _{X_i \phi X_j}Z=-\nabla _{[X_i, \phi X_j]}Z
$$
or equivalently
$$
\phi [X_i, \phi X_j] +\phi \h [X_i, \phi X_j]=0 \, .
$$
Applying $\phi$ and recalling that  $\h [X_i, \phi X_j]$ is in the kernels of $\alpha_1$ and $\alpha_2$, we get
$$
-[X_i, \phi X_j] +\mu Z_1 +\nu Z_2= \h [X_i, \phi X_j] \, ,
$$
for some functions $\mu ,\nu$.

Taking the scalar product with $X_l$ and using the symmetry of $\h$, we obtain
$$
-g([X_i, \phi X_j], X_l) =g( \h [X_i, \phi X_j] , X_l)=g( [X_i, \phi X_j] ,\h X_l)=g(  [X_i, \phi X_j] , X_l) \, ,
$$
and finally
\begin{equation}\label{eq7-th:flat-metrics}
g(  [X_i, \phi X_j] , X_l)=0\, .
\end{equation}
Similarly we have
\begin{equation}\label{eq8-th:flat-metrics}
\begin{split}
&g(  [X_i, \phi X_j] , Y_s)=g(  [Y_r, \phi Y_s] , X_l)=0\, , \\
&g(  [Y_r, \phi X_j] , X_i)=g(  [Y_r, \phi X_j] , Y_s)=0\, , \\
&g(  [\phi  Y_r, X_i] , X_l)=g(  [Y_r, \phi Y_s] , Y_t)=g(  [X_i, \phi Y_r] , Y_s)=0\, .
\end{split}
\end{equation}

\medskip

\item[c)] Now we want to show that $\nabla_{X_i}X_j \in [+1]\oplus \bR Z_1 \oplus \bR Z_2$. In fact, using Corollary \ref{c:rel-AB} and \eqref{eq2-th:flat-metrics} we have
\begin{equation*}
\begin{split}
0&=g \left ((\nabla_{\phi X_i}\phi)X_j,\phi X_l \right )=-g \left (\phi(\nabla_{\phi X_i}X_j),\phi X_l \right )\\
&=-g \left (\nabla_{\phi X_i}X_j, X_l \right) \, ,
\end{split}
\end{equation*}
where we have used the fact that $X_l\in \ker \alpha_1 \cap \ker \alpha_2$.
Since $g(X_j, \phi X_l)=0$, we obtain:
\begin{equation*}
\begin{split}
g \left (\nabla_{ X_i}X_j,\phi X_l \right )=&-g \left (X_j,\nabla_{ X_i}\phi X_l \right )\\
=&-g \left (X_j,\nabla_{ \phi X_l}X_i+[X_i , \phi X_l] \right )\\
=&-g(X_j, [X_i , \phi X_l])\\
=&0\, ,
\end{split}
\end{equation*}
where in the last equality we have used \eqref{eq7-th:flat-metrics}.

Similarly we obtain $g(\nabla_{ \phi Y_r}X_i,X_j)=0$ and in turn $g(\nabla_{ X_i}X_j,\phi Y_r)=0$. Therefore $\nabla_{X_i}X_j \in [+1]\oplus \bR Z_1 \oplus \bR Z_2$ as desired.

\medskip

\item[d)] We want to show that $[+1]$ is integrable. For $X,Y,W\in [+1]$ we have:
\begin{equation}
\begin{split}
0=&R_{XY}Z=\nabla_X (-2 \phi Y)+\nabla_Y (2 \phi X)+\phi [X,Y]+\phi \h [X,Y]\\
=&-2(\nabla_X \phi)Y+2(\nabla_Y \phi)X-\phi[X,Y]-\h \phi [X,Y] \, .
\end{split}
\end{equation}
Taking the inner product with $W$ and noting that $\phi W$ is in $[-1]$, we obtain
\begin{equation}
0=g(\phi W, [X,Y]) \, ,
\end{equation}
which means that $[X,Y]$ is orthogonal to $[-1]$.
Now let $Y\in [+1]
$ and $X\in [-1]$. By the integrability of $[-1]$, we have:
\begin{equation}
\begin{split}
0=& d\alpha_i (X, \phi Y)\\
=&-d\alpha_i (\phi X, Y)\\
=&\frac{1}{2} \alpha_i ([\phi X, Y]) \, .
\end{split}
\end{equation}
Then the bracket of two vector fields in $[+1]$ is also orthogonal to the vertical subbundle and this in turn implies that $[+1]$ is integrable.

\medskip
\item[e)] Now we prove that  $[+1]$ is totally geodesic. To show this, we will prove that, for every $Y\in[+1]$ and every $X$ orthogonal to $[+1]$, $\nabla _Y X$ is orthogonal to $[+1]$.

Let us start with $Y\in[+1]$ and $X\in [-1]$ and compute
\begin{equation}\label{eq:[+1]-tot-geod}
\begin{split}
0=&R_{XY}Z=-2\nabla_X \phi Y+\phi [X,Y]+\phi \h [X,Y]\\
=&-2(\nabla_X \phi)Y-\phi \nabla_XY-\phi \nabla_Y X -\h \phi \nabla_X Y+\h \phi \nabla_Y X \, .
\end{split}
\end{equation}
Taking the scalar product with $W\in [-1]$, we obtain
\begin{equation}
\begin{split}
0=&-2g((\nabla_X \phi)Y,W)-g(\phi \nabla_XY,W)-g(\phi \nabla_Y X ,W)-g(\phi \nabla_X Y,\h W)+g(\phi \nabla_Y X,\h W) \\
=&-2g((\nabla_X \phi)Y,W)-2g(\phi \nabla_Y X ,W) \, ,
\end{split}
\end{equation}
which implies 
$$
g((\nabla_X \phi)Y,W)=-g(\phi\nabla_Y X,W)=g(\nabla_YX,\phi W) \, .
$$
Since $g((\nabla_X \phi)Y,W)=0$ by Corollary \ref{c:rel-AB}, we have $\nabla _Y X\in [-1]\oplus \mathcal V$.

Next observe that for $X,Y$ horizontal, we have:
\begin{equation}
\begin{split}
2g((\nabla_X \phi)Z_i,Y)=&g(N^1(Z_i, Y), \phi X)-2d\alpha_i(\phi Y,X)\\
=&g(\phi^2[Z_i, Y]-\phi [Z_i, \phi Y], \phi X)-2d\alpha_i(\phi Y,X)\\
=&g(-\phi(\LL _{Z_i}\phi )Y, \phi X)-2d\alpha_i(\phi Y,X)\\
=&-g(( \LL _{Z_i}\phi)Y,  X)-2d\alpha_i(\phi Y,X) \, ,
\end{split}
\end{equation}
then $g((\nabla_X \phi)Z_i,Y)$ is symmetric in $X$ and $Y$.

Now for $Y\in[+1]$, $X\in [-1]$, taking the scalar product of \eqref{eq:[+1]-tot-geod} with $Z_i$ ($i=1,2$), we get:
\begin{equation}
\begin{split}
0=&g((\nabla_X \phi)Y, Z_i)=g((\nabla_Y \phi)X, Z_i)\\
=&-g((\nabla_Y \phi) Z_i, X)=g(\phi \nabla_Y {Z_i}, X)\\
=&-g((\nabla_Y {Z_i}, \phi X) \, ,
\end{split}
\end{equation}
which means that $\nabla_Y {Z_i}$ is orthogonal to $[+1]$.

This implies that $[+1]$ is totally geodesic and then, since 
$[-1]\oplus \bR Z_1 \oplus \bR Z_2$ is integrable with totally geodesic leaves,
the manifold splits as a local Riemannian product.

\medskip
\item[f)] 
Using equation \eqref{eq:cor1:gen-metric} and  the integrability of
$ [-1] \oplus \bR Z_1 \oplus \bR Z_2$ with totally geodesic leaves we have
\begin{equation}
\begin{split}
2g((\nabla_{X_i}\phi)X_j, Z_2)=&g( [Z_2, \phi X_j] -\phi [Z_2, X_j], X_i)+2 g(X_j, X_i)\\
=&g(  [Z_2, X_j],\phi X_i)+2 g(X_j, X_i)\\
=&g(  \nabla_{Z_2} X_j -\nabla_{X_j} Z_2, \phi X_i)+2 g(X_j, X_i)\\
=&-g( X_j,\nabla_{Z_2}  \phi X_i) +g(Z_2,\nabla_{X_j} \phi  X_i)+2 g(X_j, X_i)\\
=&g(Z_2,\nabla_{X_i} \phi  X_j)+2 g(X_j, X_i)\\
=&g(Z_2,(\nabla_{X_i} \phi )  X_j)+2 g(X_j, X_i)\\
\end{split}
\end{equation}
which gives 
\begin{equation}\label{eq10-th:flat-metrics}
g \left ((\nabla_{X_i}\phi)X_j, Z_2 \right )=2g(X_i , X_j)  \, .
\end{equation}
Using \eqref{eq10-th:flat-metrics} we obtain:
\begin{equation}\label{eq10-1-th:flat-metrics}
\begin{split}
g \left ((\nabla_{X_i}\phi)X_j, Z_1 \right)=&g \left (\nabla_{X_i}\phi X_j, Z_1 \right) \\
=&-g \left (\phi X_j, \nabla_{X_i}Z_1 \right) \\
=&-g \left (\phi X_j, \nabla_{X_i}Z \right)+g \left (\phi X_j, \nabla_{X_i}Z_2 \right) \\
=&-g \left (\phi X_j, -2 \phi X_i \right) - g \left (\nabla_{X_i}\phi X_j, Z_2 \right) \\
=&2 g \left (X_j, X_i \right) - g \left ((\nabla_{X_i}\phi )X_j, Z_2 \right) \\
=&0
\end{split}
\end{equation}
Equations \eqref{eq10-th:flat-metrics}, \eqref{eq10-1-th:flat-metrics} and Corollary \ref{c:rel-AB} then give 
 for every $X,Y \in [+1]_2$ 
\begin{equation}\label{eq11-th:flat-metrics}
(\nabla_{X}\phi)Y=2g(X, Y) Z_2 \, .
\end{equation}

Moreover \eqref{eq6-th:flat-metrics} implies $[X_i,X_j]=0$ and then $[+1]_2$ is integrable. With a similar argument we see that $[+1]_1$ is also integrable.

\medskip

\item[g)] Set $\dot h_2=\frac{1}{2}\LL _{Z_2} \phi_2$, let $L$ be a leaf of $T\FF_1$ considered as a submanifold and $X$ tangent to $L$. Let $\dot \nabla$ be the induced connection, $\sigma$ the second fundamental form, $A_{Z_1}$ the Weingarten operator in the direction $Z_1$ and $\nabla'$ the connection in the normal bundle. Then we have:
$$
\nabla_X Z= -\phi X -\phi \h X=-\phi_2 X -\phi_2 \h_2 X
$$
$$
\nabla_X Z=\nabla_X Z_1+\nabla_X Z_2=-A_{Z_1}X +\nabla'_XZ_1 +\dot \nabla _X Z_2 +\sigma(X,Z_2)\, .
$$
Comparing the previous equations we obtain
$$
A_{Z_1}X=\phi_2 (\h_2 X -\dot \h_2 X)=\phi (\h X -\dot \h_2 X)
$$
$$
 \nabla' _X Z_1=-\sigma(X,Z_2) \, .
$$
Moreover it is clear that we have $\sigma(Z_2,Z_2)=0$ and $A_{Z_1} Z_2=0$.

We want to prove that $\h_2 =\dot \h_2$ and this is equivalent to proving the vanishing of $A_{Z_1}$. To do this we will prove the vanishing of $g(A_{Z_1 }X, Y)$ for $X,Y$ elements of the eigenbasis of $\h$ constructed before.

First observe that 
$$
g(A_{Z_1 }X, Z_2)=g(\phi (\h X -\dot \h_2 X), Z_2)=0 \, ,
$$
and then we have to prove our statement for $X,Y$ horizontal. Also note that, by direct calculation, we have:
$$
g(A_{Z_1 }\phi X, Y)=g(A_{Z_1 }X, \phi Y) \, .
$$
Then it remains to prove that $g(A_{Z_1 }X_i, X_j)$ and $g(A_{Z_1 }\phi X_i, X_j)$ both vanish.

For the first one:
\begin{equation}
\begin{split}
g(A_{Z_1 }X_i, X_j)=&g(\phi (\h X_i -\dot \h_2 X_i), X_j)\\
=&g(\phi (X_i -\dot \h_2 X_i), X_j)\\
=&g(\dot \h_2 X_i,\phi X_j)\\
=&\frac{1}{2} (g([Z_2,\phi_2 X_i]-\phi_2 [Z_2, X_i],\phi_2 X_j))\\
=&\frac{1}{2} (g(\nabla_{Z_2}\phi_2 X_i-\nabla_{\phi_2 X_i} Z_2, \phi_2 X_j)-g(\phi_2 [Z_2, X_i],\phi_2 X_j))\\
=&\frac{1}{2} (g(\nabla_{Z_2}\phi_2 X_i, \phi_2 X_j)-g([Z_2, X_i],X_j))\\
=&0 \, ,
\end{split}
\end{equation}
because $[Z_2,X_i]$ is in $[-1]_2+\bR Z_2$ and where we have used
\eqref{eq1.5-th:flat-metrics} and \eqref{eq2-th:flat-metrics}.

Now we have:
\begin{equation}
\begin{split}
g(A_{Z_1 }X_i, \phi X_j)=&g(\phi_2 (\h_2 X_i -\dot \h_2 X_i), \phi_2 X_j)\\
=&g(X_i -\dot \h_2 X_i, X_j)\\
=&g(X_i , X_j)-g(\dot \h_2 X_i,X_j)\\
=&g(X_i , X_j)-\frac{1}{2} g([Z_2,\phi_2 X_i]-\phi_2 [Z_2 ,X_i],X_j)\\
=&g(X_i , X_j)-\frac{1}{2} g([Z_2, X_i],\phi_2 X_j))\\
=&g(X_i , X_j)-\frac{1}{2} g(\nabla_{Z_2} X_i -\nabla_{X_i}Z_2,\phi_2 X_j))\\
=&g(X_i , X_j)-\frac{1}{2}( -g(X_i,\nabla_{Z_2} \phi_2 X_j) +g(Z_2,(\nabla_{X_i}\phi_2) X_j))\\
=&0\, ,
\end{split}
\end{equation}
where at the end we used \eqref{eq1.5-th:flat-metrics} and \eqref{eq11-th:flat-metrics}.

\medskip

\item[h)] To prove that  $[+1]_2$ is also totally geodesic, let us consider the operators $H_i=\frac{1}{2} \LL_{Z_i}\phi$. Each $H_i$ is symmetric by Theorem 3.4.  From our observation above that 
$\h_2 =\dot \h_2$ a simple direct computation shows that we have
$$H_2 X_i=X_i\;{\rm and \;similarly}\;  H_1 Y_r=Y_r.$$ 
This implies that  $H_1 X_i$ has no $Y_r$ or $\phi Y_r$ component.  Thus since  $A_{Z_1 }$ vanishes we have
 $$H_1 X_i=0\;{\rm and \;similarly}\; H_2 Y_i=0.$$

Applying \eqref{eq11-th:flat-metrics}, we obtain
$$
2g(X_i, X_j)Z_2= \nabla_{X_i} \phi X_j-\phi \nabla _{X_i}X_j \, ,
$$
which, differentiating along $Y_r$, gives:
\begin{equation}
2(Y_r g(X_i, X_j))Z_2+2 g(X_i, X_j)\nabla_{Y_r}Z_2= \nabla_{Y_r}\nabla_{X_i} \phi X_j-(\nabla_{Y_r}\phi) \nabla _{X_i}X_j -\phi \nabla_{Y_r}\nabla _{X_i}X_j \, .
\end{equation}
Taking the scalar product with $Z_2$, we get
\begin{equation}\label{eq:tot-geod-[+1]1-1}
2Y_r g(X_i, X_j)= g(\nabla_{Y_r}\nabla_{X_i} \phi X_j,Z_2)-g((\nabla_{Y_r}\phi) \nabla _{X_i}X_j ,Z_2) \, .
\end{equation}

Now, applying \eqref{eq:cor1:gen-metric}, we obtain
\begin{equation}
\begin{split}
-g((\nabla_{Y_r}\phi) \nabla _{X_i}X_j ,Z_2) =&g(\nabla _{X_i}X_j ,(\nabla_{Y_r}\phi)Z_2) \\
=&-g(\phi H_2 \nabla _{X_i}X_j ,\phi Y_r) -d\alpha_2(\phi \nabla _{X_i}X_j, Y_r)\\
=&-g(H_2 \nabla _{X_i}X_j , Y_r)\\
=&-g(\nabla _{X_i}X_j , H_2Y_r)\\
=&0\, .
\end{split}
\end{equation}
Applying \eqref{eq:cor1:gen-metric} again we first note that by computing as we have been doing, the above properties of $H_1$ and $H_2$ yield
$$g((\nabla_{Y_r}\phi)X_j,Z_1)=g((\nabla_{Y_r}\phi)X_j,Z_2)=0.$$
Corollary \ref{c:rel-AB} then gives $(\nabla_{Y_r}\phi)X_j=0$ and therefore by  \eqref{eq:cor1:gen-metric} 
\begin{equation}
\begin{split}
g(\nabla_{X_i}\nabla_ {Y_r}\phi X_j,Z_2)=&g((\nabla_{X_i}\phi) \nabla_ {Y_r}X_j,Z_2)\\
=&+g(H_2\nabla_ {Y_r}X_j,X_i)+d\alpha_2(\phi \nabla_ {Y_r}X_j, X_i)\\
=&2g(\nabla_ {Y_r}X_j,X_i)\, .
\end{split}
\end{equation}

Equation \eqref{eq:tot-geod-[+1]1-1} then becomes
\begin{equation}\label{eq:tot-geod-[+1]1-2}
\begin{split}
2Y_r g(X_i, X_j)=&g(\nabla_{Y_r}\nabla_{X_i} \phi X_j,Z_2)\\
=&g(R_{Y_rX_i}\phi X_j,Z_2)+g(\nabla_{X_i}\nabla_{Y_r} \phi X_j,Z_2)\\
=&2 g(\nabla_ {Y_r}X_j,X_i)
\end{split}
\end{equation}
since $[Y_r,X_j]=0$ by \eqref{eq6-th:flat-metrics}.
 But, by the compatibility condition of the Levi-Civita connection we also have
$$
Y_r g(X_i, X_j)=g(\nabla_{Y_r} X_i, X_j)+g(X_i, \nabla_{Y_r} X_j) \, .
$$
Comparing this with \eqref{eq:tot-geod-[+1]1-2}, gives
$$
g(\nabla_ {Y_r}X_i
,X_j)=0 \, ,
$$
and in turn, again noting $[Y_r,X_j]=0$, 
$$g(Y_r,\nabla_{X_j} X_i)=0.$$
Since $[+1]$ is totally geodesic, we also have $g(\nabla_{X_j} X_i,\phi Y_r)=0$,
therefore  $[+1]_2$ is totally geodesic.

\medskip

\item[i)] We rewrite \eqref{eq11-th:flat-metrics} as follows
\begin{equation}\label{eq14-th:flat-metrics}
2g(X_i, X_j) Z_2=(\nabla_{ X_i}\phi)X_j=\nabla_{ X_i} \phi X_j - \phi \nabla_{ X_i} X_j \, ,
\end{equation}
and we want to apply $\nabla _{X_l}$ to \eqref{eq14-th:flat-metrics}. We firstly need the following calculations:
\begin{equation}
\begin{split}
2g(X_i , X_j)=&g \left ((\nabla_{X_i}\phi)X_j, Z_2 \right )\\
=&- g \left (X_j,(\nabla_{X_i}\phi) Z_2 \right )\\
=&- g \left (X_j,-\phi(\nabla_{X_i} Z_2) \right ) \\
=&- g \left (\phi X_j,\nabla_{X_i} Z_2 \right ) \, .
\end{split}
\end{equation}
and therefore 
\begin{equation*}
g(\nabla_{X_l}Z_2, \phi X_m)=-2g(X_l, X_m)\, .
\end{equation*}

Applying $\nabla _{X_l}$ to  \eqref{eq14-th:flat-metrics} we obtain
\begin{equation*}
\begin{split}
2(X_l g(X_i, X_j))Z_2+2 g(X_i, X_j)\nabla_{ X_l}Z_2=\\
\nabla_{ X_l}\nabla_{X_i}\phi X_j-(\nabla_{X_l} \phi)(\nabla_{ X_i}X_j)-\phi( \nabla_{X_l}\nabla_{ X_i}X_j) \, .
\end{split}
\end{equation*}
Taking the scalar product with $\phi X_m$ we get
\begin{equation}\label{eq18-th:flat-metrics}
-4g(X_i,X_j) g(X_l,X_m)=g(\nabla_{ X_l}\nabla_{X_i}\phi X_j-\phi( \nabla_{X_l}\nabla_{ X_i}X_j), \phi X_m) \, ,
\end{equation}
since $g((\nabla_{X_l} \phi)(\nabla_{ X_i}X_j), \phi X_m)=0$ by Corollary \ref{c:rel-AB}.

\medskip

\item[j)] Note that from Theorem \ref{th:ass-metrics} we know that the integral submanifolds of $\mathcal{V}$ are totally geodesic. Moreover we proved that $[+1]_1$ and $[+1]_2$ are integrable and totally geodesic. We also know that  $[-1]_1$ and $[-1]_2$ are integrable, totally geodesic and flat by \eqref{eq2-th:flat-metrics} and \eqref{eq3-th:flat-metrics}. Then the manifold $M$ splits locally as a Riemannian product of integrable submanifolds of $\mathcal{V}$, $[\pm 1]_1$ and $[\pm 1]_2$.  To conclude the proof we exhibit the curvature. To do this interchange the role of $X_l$ and $X_i$ in \eqref{eq18-th:flat-metrics} and take the difference. Then we have
\begin{equation*}
\begin{split}
-4g(X_i,X_j) g(X_l,X_m)+4g(X_l,X_j) g(X_i,X_m)=
g(R_{X_lX_i}\phi X_j,\phi X_m)-g(R_{X_lX_i}X_j, X_m).
\end{split}
\end{equation*}
The first term on the right vanishes by the local Riemannian product structure and the
 second term then gives the desired value of the curvature. If $h$ or $k$ is equal to 1, the corresponding $[+1]_i$ subbundle is 1-dimensional and the leaves of the characteristic foliation are flat.
\end{noindenum}

\end{proof}
\begin{corollary}
Let $M$ be a $(2h+2k+2)$-dimensional manifold endowed with a metric contact pair $(\alpha_1, \alpha_2, \phi, g)$ of type $(h,k)$ (with $h\geq 1$) and decomposable $\phi$, and such that the curvature of  $g$ vanishes on the vertical subbundle. If $M$ is complete, then its Riemannian universal covering is isometric to $\bE^{h+1}\times \bS^h(4)\times \bE^{k+1}\times \bS^k(4)$ if $k\geq 1$ or $\bE^{h+1}\times \bS^h(4)\times \bE^{1}$ if $k=0$.  It is to be understood that when $h$ (or $k$) is equal to 1, the $\bS^h(4)$ factor will just contribute another line to the Euclidean factor.
\begin{proof}
The Riemannian universal covering $\tilde M$ is locally isometric to $M$ and then by Theorem \ref{th:flat-metrics} is locally isometric to $\bE^{h+1}\times \bS^h(4)\times \bE^{k+1}\times \bS^k(4)$ if $k\geq 1$ or to $\bE^{h+1}\times \bS^h(4)\times \bE^{1}$ if $k=0$. Then one concludes by applying the de Rham Decomposition Theorem.
\end{proof}

\end{corollary}

\begin{corollary}
Let $M$ be a $(2h+2k+2)$-dimensional manifold endowed with a metric contact pair $(\alpha_1, \alpha_2, \phi, g)$ of type $(h,k)$ (with $h\geq 1$) and decomposable $\phi$. If $g$ is flat then $h,k \leq 1$.
\end{corollary}

\bibliographystyle{amsplain}

\begin{thebibliography}{999}


\bibitem{Abe}
K.~Abe, \emph{On a class of Hermitian manifolds},  Invent. Math.  \textbf{51}  (1979),
103--121.


\bibitem{Bande1}
G.~Bande, \textsl{Formes de contact g{\'e}n{\'e}ralis{\'e},
couples de contact et couples contacto-symplectiques}, Th{\`e}se
de Doctorat, Universit{\'e} de Haute Alsace, Mulhouse, 2000.

\bibitem{Bande2}
G.~Bande, \emph{Couples contacto-symplectiques},
Trans.~Amer.~Math.~Soc. \textbf{355} (2003),
1699--1711.

\bibitem{BGK} G.~Bande, P.~Ghiggini and D.~Kotschick, \emph{Stability theorems
for symplectic and contact pairs}, Int. Math. Res. Not. \textbf{68} (2004), 3673--3688.

\bibitem{BH} G.~Bande and A.~Hadjar, \emph{Contact Pairs}, Tohoku Math. J. \textbf{57} (2005),
247--260.

\bibitem{BH2}
G.~Bande and A.~Hadjar, \emph{Contact pair structures and associated
metrics}, Differential Geometry - Proceedings of the 8th International Colloquium, World Sci. Publ. (2009), 266--275.

\bibitem{BH3}
G.~Bande and A.~Hadjar, \emph{On normal contact pairs}, Internat.~J.~Math.~\textbf{21} (2010),
737--754.

\bibitem{BH4}
G.~Bande and A.~Hadjar, \emph{On the characteristic foliations of metric contact pairs}, Harmonic maps and differential geometry, 255--259, Contemp. Math., \textbf{542}, Amer. Math. Soc., Providence, RI, 2011.



\bibitem{BK} G.~Bande and D.~Kotschick, \emph{The Geometry of Symplectic pairs},
Trans.~Amer.~Math.~Soc.~\textbf{358} (2006),
1643--1655.

\bibitem{BK2} G.~Bande and D.~Kotschick, \emph{Contact pairs and locally conformally symplectic structures}, Harmonic maps and differential geometry, 85--98, Contemp. Math., \textbf{542}, Amer. Math. Soc., Providence, RI, 2011..




\bibitem{Blair3}
D.~E.~Blair, \emph{On the non-existence of flat contact metric structures},  Tohoku Math. J. \textbf{28} (1976),
373--379.

\bibitem{Blair4} 
D.~E.~Blair, \emph{Two remarks on contact metric structures},  Tohoku Math. J. \textbf{29} (1977), 319--324.



\bibitem{Blairbook}
D.~E.~Blair, \emph{Riemannian geometry of contact and symplectic
manifolds}, Progress in Mathematics, 2nd Ed., vol. 203, Birkh\"auser, 2010.

\bibitem{Blair2}
D.~E.~Blair, G.~D.~Ludden and K.~Yano,  \emph{Geometry of complex manifolds similar to the Calabi-Eckmann manifolds},  J. Differential Geom.  \textbf{9}  (1974), 263--274.


\bibitem{Bor}
A.~Borel, \emph{Compact Clifford--Klein forms of symmetric spaces},
Topology {\bf 2} (1963), 111--122.


\bibitem{Ornea}
S.~Dragomir and L.~Ornea, \emph{Locally conformal K\"ahler geometry},
Progress in Mathematics, vol. 155, Birkh\"auser, 1998.

\bibitem{Geiges2}
H.~Geiges, \emph{Normal contact structures on $3$-manifolds},
Tohoku Math. J. {\bf49} (1997),
415--422.

\bibitem{HOT}
Y.~Hatakeyama, Y.~Ogawa and S.~Tanno, \emph{Some properties of manifolds with contact metric structure}, Tohoku Math. J. \textbf{15} (1963), 42--48.

\bibitem{Vaisman}
I.~Vaisman, \emph{Some curvature properties of locally conformal K\"ahler manifolds},
Trans.~Amer.~Math.~Soc. \textbf{259} (1980), 439--447.


\end{thebibliography}

\end{document}